\documentclass[12pt]{article}
\usepackage[english]{babel}
\usepackage{color}
\usepackage[dvipsnames]{xcolor}
\usepackage{amsmath}
\usepackage{physics}
\usepackage{xr-hyper}
\usepackage[hidelinks]{hyperref}
\usepackage{diagbox}
\usepackage{multirow}
\usepackage{amsfonts}
\usepackage{amssymb}
\usepackage{latexsym}
\usepackage{graphicx}
\usepackage{adjustbox}
\usepackage{float}
\usepackage{xr}
\usepackage{caption}
\usepackage{enumitem}
\usepackage{subfigure}
\usepackage{seqsplit}
\usepackage[titletoc]{appendix}
\usepackage[normalem]{ulem}
\usepackage[margin=0.75in]{geometry}
\usepackage{ulem}
\usepackage{cancel}
\usepackage{xcolor}

\def\theequation{\arabic{section}.\arabic{equation}}

\usepackage{mathtools}
\usepackage[nodisplayskipstretch]{setspace} 

\newcommand{\pa}{\partial}

\newcommand{\ta}{\tau}

\newcommand{\rar}{\rightarrow}

\newcommand{\ii}{\mathrm{i}}
\newcommand{\ee}{\mathrm{e}}
\newcommand{\D}{\mathrm{d}}

\newtheorem{thm}{Theorem}[section]
\newtheorem{lemma}[thm]{Lemma}

\newtheorem{assumption}[thm]{Assumption}

\usepackage[nottoc]{tocbibind}

\begin{document}

\title{Fourth-order compact exponential splittings for unbounded operators}
	\author{ Juan Carlos Del Valle\footnote{Faculty of Mathematics, Physics, and Informatics, University of Gda\'nsk, Gda\'nsk, Poland.}  \qquad Arieh Iserles \footnote{DAMTP, University of Cambridge, Cambridge, UK.}
\qquad Karolina Kropielnicka\footnote{Institute of Mathematics, Polish Academy of Sciences, Warsaw, Poland.}}
	\date{\today}
	
		\maketitle
	\begin{abstract}
	We present a derivation and error bound for the family of fourth order splittings, originally introduced by Chin and Chen, where one of the operators is unbounded and the second one bounded but time dependent, and which are dependent on a parameter. We first express the error by an iterated application of the Duhamel principle, followed by quadratures of Birkhoff-Hermite type of the underlying multivariate integrals. This leads to error estimates and bounds, derived using Peano/Sard kernels and direct estimates of the leading error term. Our analysis demonstrates that, although no single value of the parameter can minimise simultaneously all error components, an excellent compromise is the cubic Gauss--Legendre point $1/2-\sqrt{15}/10$.
	\end{abstract}

\section{Introduction}

Exponential splittings are a major tool in the solution of time-dependent partial differential equations, in particular, in situations when it is important to preserve qualitative features of the underlying solution \cite{Blanes,McLachlan}. Order analysis of splitting methods are straightforward when  differential operators are all bounded, all it takes is a Taylor expansion. However, in the presence of unbounded operators Taylor expansions no longer converge and the determination of order becomes considerably more challenging. This has been addressed in a number of  papers  using functional-analytic machinery \cite{FamilyStrang,hansen09esu,Lubich}. 

A prominent role in the practice of splitting methods, in particular for the linear Schr\"odinger equation
\begin{equation}
  \label{Shreddy}
  \mathrm{i}\frac{\partial u}{\partial t} =-\frac12 \Delta u +V(x,t)u,\qquad u(0,t)=u_0,\quad x\in\Omega\subseteq\mathbb{R}^d,
\end{equation}
given with suitable boundary conditions, is played by {\em compact splittings,\/} also known as {\em force gradient splitting methods},  introduced by Chen and Chin in \cite{chin02gsa} and Chin and Anisimov in \cite{Chin}.
In particular, they allow to attain order~4 by very efficient (indeed, compact) schemes  that guarantee the celebrated \textit{mass preservation}. However, their method of proof requires that the underlying operators are all bounded and this is not the case once $\Omega=\mathbb{R}^d$, a state of affairs of major importance in quantum mechanics. 
     
Let $A$ and $B(t)$, $t\geq0$, be linear operators and suppose that $A:D(A)\subset X\rar X$ generates a strongly continuous semigroup. While $A$ might be unbounded, we stipulate that $B(t)$ is bounded for every $t\geq0$. It is instructive to bear in mind the example (\ref{Shreddy}), whereby $A=\frac12 \mathrm{i}\Delta$, while $B(t)$ is the bounded multiplication operator $-\mathrm{i}V(x,t)$. Thus, our concern is with the general linear-evolutionary equation
\begin{equation}
  u'(t) =[A+B(t)]u(t) ,\qquad u(0)=u_0\qquad  t\in[0,T],\quad x\in\Omega,
  \label{main}
\end{equation}
where we assume (i) $[B(t),B'(t)] = 0$ and (ii) $[B(t),[A,B(t)]]$ is a multiplication operator.
These  assumptions are  fulfilled  in the case of the  linear Schr\"odinger equation (\ref{Shreddy}). In particular, the commutator \hbox{$[B(t),[A,B(t)]]=-\ii \sum_k| \pa_{x_k} V(x,t)|^2$} represents the squared norm of the {\em force}.

Note that, being dependent on $t$, the operator $B(t)$ does not generate a semigroup. However, $B(\tilde{t})$  trivially generates a semigroup for frozen $\tilde{t}$, and we will make use of this in the sequel.


One of the proposed fourth-order algorithms is the $\tau$-dependent {\em family} $\mathcal{T}_{ACB}^{(4)}$, presented in \cite{chin02gsa}, see equation (22) therein. It has the form
\begin{equation}
\label{TACB}
u(h) \approx \mathcal{T}_{ACB}^{(4)}u_0:=\ee^{\tau h A}\,\ee^{p hB((1-\tau)h)}\,\ee^{(\frac12-\tau)hA}\,\ee^{h\tilde{B}(h/2)}\,\ee^{(\frac12-\tau)hA}\,\ee^{ph B(\tau h)}\,\ee^{\tau hA}u_0,
\end{equation}	
where $\tau \in[0,\frac12)$ is a parameter,
\begin{equation}
\label{Btilde}
	\tilde{B}(t)\ =\ qB(t)+r\,h^2[B(t),[A,B(t)]],
\end{equation}
and
\begin{equation}
  \label{coeffs}
p=\frac{1}{6(1-2\tau)^2},\qquad q=1-2p=\frac{2 \left(6 \tau ^2-6 \tau +1\right)}{3 (2 \tau -1)^2},\qquad r=\frac{1}{12}\left[1-\frac{1}{1-2\tau}+\frac{1}{6(1-2\tau)^3}
 \right]\!.
 \end{equation}
Because of our assumptions, $\tilde{B}$ is a bounded operator. 
Note that by taking $\ta=0$, the outermost semigroups in (\ref{TACB}) become identities, reducing the number of exponentials from seven to five. There is no other value of $\ta\in[0,\frac12)$ with this property.

In fact, authors of \cite{chin02gsa}  not only proposed the fourth-order {\em family} $\mathcal{T}_{ACB}^{(4)}$ but also {\em family} $\mathcal{T}_{BDA}^{(4)}$, labeled in their manuscript with number (33). The main difference between them is their structure. In contrast to $\mathcal{T}_{ACB}^{(4)}$, $\mathcal{T}_{BDA}^{(4)}$  is characterized by having four exponentials involving only multiplication operators and three more comprising just of $A$. The derivation of  both families is based on so called {\em forward time derivative}, also known as {\em super operator}, see \cite{Suzuki}. Given the intuitive yet heuristic arguments they have used,  error analysis of both families cannot be followed in detail. 
The convergence is open to challenge, especially in the case of an unbounded operator $A$ (e.g.\ a Laplacian in an unbounded domain). However, numerical experiments found in the literature \cite{Baye,chin02gsa,KIERI201533} support the conclusion that they are indeed fourth-order exponential splittings.
Furthermore,  rigorous mathematical study of the super operator is missing to the best of present authors' knowledge.
Consequently, neither  convergence analysis nor error bounds of families  $\mathcal{T}_{ACB}^{(4)}$ and $\mathcal{T}_{BDA}^{(4)}$ have been ever derived in a rigorous manner.

While  progress concerning the convergence was made for time-independent $B$  in 
 \cite{KIERI201533}, the more general time-dependent case  remains outstanding. Also, the choice of optimal $\tau$, the one that minimizes the error constant of the integrator in families $\mathcal{T}_{ACB}^{(4)}$ and $\mathcal{T}_{BDA}^{(4)}$, has been never determined. 
This merits further investigation and has been a major motivation for the present work.


Under similar assumptions, a simpler family of second-order exponential splittings of the form
\begin{equation}
\label{familyII}
u(h)\approx \ee^{h(1-\tau)A}\ee^{hB(h\tau)+\frac{h^2(1-2\tau)}{2}\mathcal{C}(h\tau)}\ee^{h\tau A} u_0
\end{equation}
where $\tau\in[0,1]$ and $\mathcal{C}(t)=[B(t),A]+B'(t)$  has been recently  investigated  in \cite{FamilyStrang}. The approximate solution defined by (\ref{familyII})  stands out due to the appearance of the commutator and a time derivative. In practice,  (\ref{familyII}) serves as a minimal example of a splitting that contains a commutator. The work presented in \cite{FamilyStrang} reveals that  determination of the order and error analysis can be done simultaneously by employing elementary tools: iteration on the variation-of-constants formula (also known as the {\em Duhamel's formula\/}) and Hermite--Birkhoff quadratures to integrate all relevant terms generated by the iteration. In fact, the presence of derivatives in the quadratures is essential since ultimately generates the commutators in the splitting.

To increase accuracy beyond second order, one may consider more exponentials involved in the splitting. In the language of \cite{FamilyStrang}, this means more iterations and higher degree of accuracy in the Hermite--Birkhoff quadratures. However, as is well known \cite{sheng89slp},  there is no splitting method of order higher than two involving positive coefficients. As an alternative, one could consider  negative \cite{hansen09hos} or even imaginary coefficients. Unfortunately, splittings for Schr\"odinger equations would  not be stable in the presence of imaginary coefficients, while splittings for parabolic equations do not converge for negative propagation in time. This is the advantage of compact splittings of type $\mathcal{T}_{ACB}^{(4)}$ and $\mathcal{T}_{BDA}^{(4)}$: both positive coefficients and a reduced number of exponentials.

In this paper, we will present the derivation and error analysis of $\mathcal{T}_{ACB}^{(4)}$,  while  family $\mathcal{T}_{BDA}^{(4)}$ can be tackled in a similar way.   The starting point is  iterating the variation-of-constants formula. The nontrivial part, a  reconstruction of the semigroup ${\rm e}^{qhB(t)+rh^2[B(t),[A,B(t)]]}$, leads to  quadratures featuring first order derivatives. This can be done, however, only once we establish a relation between the derivatives of the integrands and the commutators. We  also discuss the influence of parameter $\tau$ on the accuracy, striving to minimise the error constant of the method. 


The present study is grounded on the following two assumptions
\begin{assumption}\label{assum:main_assum}
	We assume that $A$ is a densely defined and closed, linear (possibly unbounded) operator $A:D(A)\subset X \rightarrow X$ that  generates of a strongly continuous $C_0$-semigroup ${\rm e}^{tA}$ on $X$, and that, for each fixed $s\in [0,h]$, $B(s)$ is a bounded, linear operator acting on $X$, that is $\forall_{s\in[0,h]}B(s)\in\mathfrak{B}(X)$, where $\mathfrak{B}(X)$ is the space of bounded linear operators defined on $X$. Moreover we assume that $B\in C^4([0,h],\mathfrak{B}(X))$ and seek solutions $u\in C^1([0,h],X)$.
\end{assumption}

\begin{assumption}\label{assum:commutators}
	We assume that 
	\begin{align*}
		(i)&\quad [B(t),B'(t)] = 0,\\
		(ii)&\quad [B(t),[A,B(t)]]\ {\rm is\ a\ multiplication\ operator}.
	\end{align*}
\end{assumption}

The paper is organized as follows. In Section \ref{section:derivation} we present the derivation of the method $\mathcal{T}_{ACB}^{(4)}$ while in Section \ref{section:ErrorTerms} we focus on the error analysis, discuss the importance of $\tau$ and formulate the theorem on the convergence of the method. 
Numerical example is presented in Section \ref{section:numerical_example} and finally conclusions are presented in Section \ref{section:conclusions}. Detailed (and fairly dense) algebra is relegated to appendices.

\setcounter{equation}{0}
\section{The derivation of the method}\label{section:derivation}

For transparency of the presentation, this section is focused just on the derivation of the family $\mathcal{T}_{ACB}^{(4)}$, while the study of  relevant error terms is postponed to Section \ref{section:ErrorTerms}.  We commence by recalling that, being bounded, the operator $B(s)$ {\em with the variable $s$ frozen\/} generates a strongly continuous semigroup. In particular, 
\begin{equation}\label{semigroup_B_expansions}
  \mathrm{e}^{chB(s)}=\sum_{n=0}^\infty \frac{(ch)^n}{n!} B^n(s)
\end{equation}
and the series converges. Because of our assumption that $[B(t),[A,B(t)]]$ is a multiplication operator, the same applies to $\tilde{B}$ (which has been given in (\ref{Btilde})) and
\begin{eqnarray}\label{semigroup_B_tilde_expansions}
  \ee^{h\tilde{B}(\frac{h}{2})}&\!\!\!=\!\!\!&\sum_{n=0}^\infty \frac{h^n}{n!} \tilde{B}^n(\tfrac{h}{2} )\\ \nonumber
  &\!\!\!=\!\!\!&I+hqB(\tfrac{h}{2})+\tfrac12h^2 q^2 B^2(\tfrac{h}{2})+h^3 \left\{r[B(\tfrac{h}{2}),[A,B(\tfrac{h}{2})]]+\tfrac16 q^3 B^3(\tfrac{h}{2})\right\}\\ \nonumber
  &\!\!\!\!\!\!&\mbox{}+h^4 \left\{ \tfrac12 qr B(\tfrac{h}{2})[B(\tfrac{h}{2}),[A,B(\tfrac{h}{2})]]+\tfrac12 qr [B(\tfrac{h}{2}),[A,B(\tfrac{h}{2})]] B(\tfrac{h}{2})+\tfrac{1}{24} q^4 B^4(\tfrac{h}{2})\right\}\\ \nonumber
  &\!\!\!\!\!\!&\mbox{}+\mathcal{O}(h^5)\\ \nonumber
  &\!\!\!=\!\!\!&I+hqB(\tfrac{h}{2})+ \tfrac12h^2 q^2 B^2(\tfrac{h}{2})+h^3 \left\{r[B(\tfrac12 h),[A,B(\tfrac{h}{2})]]+\tfrac16 q^3 B^3(\tfrac{h}{2})\right\}\\ \nonumber
  &\!\!\!\!\!\!&\mbox{}+h^4 \left\{ \tfrac12 qr [B^2(\tfrac{h}{2}),[A,B(\tfrac{h}{2})]]+\tfrac{1}{24} q^4 B^4(\tfrac{h}{2})\right\}+\mathcal{O}(h^5). \nonumber
\end{eqnarray}

This means that the splitting (\ref{TACB}) is approximated by
\begin{eqnarray} \label{expansion}
u(h)&\!\!\!\approx\!\!\!&\ee^{h\tau A} \times \sum_{n=0}^\infty \frac{(hp)^n}{n!} B^n((1-\tau)h) \times \ee^{h(\frac12-\tau)A} \times \sum_{n=0}^\infty \frac{h^n}{n!} \tilde{B}^n(\tfrac12 h) \times \ee^{h(\frac12-\tau)A} \\ \nonumber
  &\!\!\!\!\!\!&\mbox{}\times \sum_{n=0}^\infty \frac{(hp)^n}{n!} B^n(\tau h) \times \ee^{h\tau A}u_0,
\end{eqnarray}
which, in turn, is a sum of terms of the form
\begin{equation}
  \label{NiceTerm}
  \ee^{h\tau A} V_1 \ee^{h(\frac12-\tau)A} V_2 \ee^{h(\frac12-\tau)A} V_3 \ee^{h\tau A}u_0,
\end{equation}
where $V_1,V_3$ consist of summands of the right hand side of (\ref{semigroup_B_expansions}) while $V_2$ consists of the summands of the right hand side of (\ref{semigroup_B_tilde_expansions}).
We  are interested  in expressions which are $\mathcal{O}(h^k)$ for $k\leq 4$. The first four terms are of order $\mathcal{O}(h^0)$ and $\mathcal{O}(h)$:
\begin{eqnarray*}
  & \ee^{hA}u_0, \\
  & hp\, \ee^{h\tau A} \times B((1-\tau)h)  \times  \ee^{h(1-\tau)A}u_0, \\
  & hq\,\ee^{\frac{1}{2} h A} \times  B(\tfrac12 h)  \times \ee^{\frac{1}{2} h A}u_0,\\ 
  & hp\,\ee^{h(1-\tau)A} \times B(\tau h) \times \ee^{h\tau A}u_0.\\
\end{eqnarray*}
Additionally, we have six terms of order $\mathcal{O}(h^2)$, ten terms of order $\mathcal{O}(h^3)$ and finally fifteen terms of order $\mathcal{O}(h^4)$.  All thirty five terms are listed in  Appendix \ref{explicitterms}.

We will now show  that all these terms correspond to the expressions originating in specific Hermite--Birkhoff quadratures used to integrate a four-times iterated Duhamel's formula. 
To do so, let us consider the mild solution to (\ref{main}) provided by the  variation-of-constants formula in the form 
\begin{equation}\label{Duhamel}
u(h)={\rm e}^{hA}u_0+\int_0^h\ee^{(h-\eta_1)A}B(\eta_1)u(\eta_1)\,\text{d}\eta_1.
\end{equation}
Iterating (\ref{Duhamel}) four times, we obtain
\begin{align}
\label{iterations}
u(h)=\ee^{hA}u_0\ &+\ \int_0^{h}f_1(\eta_1)\,\text{d}\eta_1
 +\ \int_0^{h}\int_0^{\eta_1}f_2(\eta_1,\eta_2)\,\text{d}\eta_2\text{d}\eta_1\nonumber\\
&+ \int_0^{h}\int_0^{\eta_1}\int_0^{\eta_2}f_3(\eta_1,\eta_2,\eta_3)\,\text{d}\eta_3\text{d}\eta_2\text{d}\eta_1\nonumber\\
& + \int_0^{h}\int_0^{\eta_1}\int_0^{\eta_2}\int_0^{\eta_3}f_4(\eta_1,\eta_2,\eta_3,\eta_4)\,\text{d}\eta_4\text{d}\eta_3\text{d}\eta_2\text{d}\eta_1               
+
R_V ,
\end{align}
where

\begin{equation}\label{fn}
f_n(\eta_1,\ldots,\eta_n) =\ee^{hA}\,\mathcal{T}\left\{\prod_{i=1}^{n}\ee^{-\eta_iA}\,B(\eta_i)\,\ee^{\eta_iA}\right\}u_0,\ n=1,\ldots,4
\end{equation}
and

\begin{equation}
  \label{RV}
R_V = \int_0^{h}\!\!\!\!\ldots\!\!\int_0^{\eta_4}\!\ee^{(h-\eta_1)A}B(\eta_1)\ee^{(\eta_1-\eta_2)A}B(\eta_2)\ee^{(\eta_2-\eta_3)A}B(\eta_3)\ee^{(\eta_3-\eta_4)A}B(\eta_4)\ee^{(\eta_4-\eta_5)A}B(\eta_5)u(\eta_5)\text{d}\eta_5\ldots\text{d}\eta_1.
\end{equation}
Here $\mathcal{T}$ is the time-ordering operator, which enforces the correct sequence of factors in the product.
Note that the  index $n$ in $f_n$ indicates the number of insertions of the operator $B(t)$ at different times, and also the multiple $n$-variable integration over the $n$-simplex  that has to be performed in (\ref{iterations}) and (\ref{RV}).
Equations (\ref{Duhamel}) and (\ref{iterations}) define the solution in an implicit manner and may be regarded as a first step towards the construction of exponential integrators, see \cite{hochbruck_ostermann_2010}. Given that $\mathcal{T}_{ACB}^{(4)}$ will be derived from (\ref{iterations}), we can consider it simultaneously as special case of an exponential integrator and a splitting method. To do so, we regard $R_V$ as a source of local error of order $\mathcal{O}(h^5)$ and look for quadratures approximating the remaining integrals. Obviously, the choice should be  consistent with fifth-order accuracy. This means, that the quadrature for the $n$-fold integrals needs to be accurate for polynomials of degree (at least)  $5-n$, $n=1,\ldots,4$. The sum of the semigroup ${\rm e}^{hA}$ and of the four correspondent quadratures is used to 
construct $\mathcal{T}^{(4)}_{ACB}$. For this reason, we expect the quadratures to be based on such evaluations of the integrands $f_n, \ n=1,2,3,4$, which  reconstruct the expansions (\ref{expansion}) up to accuracy of $\mathcal{O}(h^5)$. More precisely we expect the evaluations of the integrands and their derivatives to correspond to terms (\ref{NiceTerm}) listed in Appendix \ref{explicitterms}. 
This suggests, for example, that the quadrature for the first integral in (\ref{iterations}) should be approximated by values of $f_1$ at the quadrature nodes 
\begin{equation}
\label{eq:1d_nodes}
  \alpha_3=(1-\tau)h,\qquad \alpha_2=\frac{1}{2}h,\qquad \alpha_1=\tau h,
\end{equation}
with quadrature weights 
\begin{displaymath}
  v_1=v_3=p,\qquad v_2=q,
\end{displaymath}
where $p$ and $q$ have been given in (\ref{coeffs}). Thus, the univariate integral is approximated as
\begin{equation}
	\int_{0}^hf_1(\eta_1)\text{d}\eta_1 =h \sum_{k=1}^3 v_{k} f_1(\alpha_k) + R_{1},
	\label{1dimensional}
\end{equation}
where  $R_{1}=\mathcal{O}(h^5)$ is the remainder term that might be easily obtained by elementary means, e.g., the Peano kernel theorem.

\subsection{A general pattern of quadrature}

The pattern (\ref{1dimensional}) does not generalise in a straightforward manner to higher dimensions, not least because, to attain order 5, it is not enough to compute the integrand just at the natural generalisation of univariate quadrature points (\ref{eq:1d_nodes}). Specifically, in dimensions 2 and 3 we need also to use derivatives for the quadratures, a procedure known as {\em Hermite--Birkhoff quadrature\/}  \cite{lorenz83bi}.

Hermite--Birkhoff quadrature guarantees the requisite order \textcolor{cyan}{5} and, in addition,  the reconstruction of the semigroup $\ee^{qB(t)+r\,t^2[B(t),[A,B(t)]]}$ (where $t$ in $B(t)$ is assumed frozen). A crucial point is that
\begin{equation}
\label{key_remark}
  \left[
  \begin{array}{cc}
         1& -1
  \end{array}
  \right]\cdot\grad f_2(\alpha_2,\alpha_2)
= \ee^{\frac{h}{2}A}[B(\alpha_2),[A,B(\alpha_2)]]\ee^{\frac{h}{2}A}u_0,
\end{equation}
as can be verified by direct calculation and by the observation that $\ee^{\frac{h}{2}A}[B'(\tfrac{h}{2}),B(\tfrac{h}{2})]\ee^{\frac{h}{2}A}=0$, a consequence of our assumption (i). Note that $[\begin{array}{cc}1&-1\end{array}]^\top$ is the normal at the midpoint $(\alpha_2,\alpha_2)$ of the diagonal of the two-dimensional simplex. Similar observation remains valid  for the triple integral, while for the quadruple integral there is no need to use derivatives to attain order 5.

\begin{figure}[htb]
  \begin{center}
    \includegraphics[width=200pt]{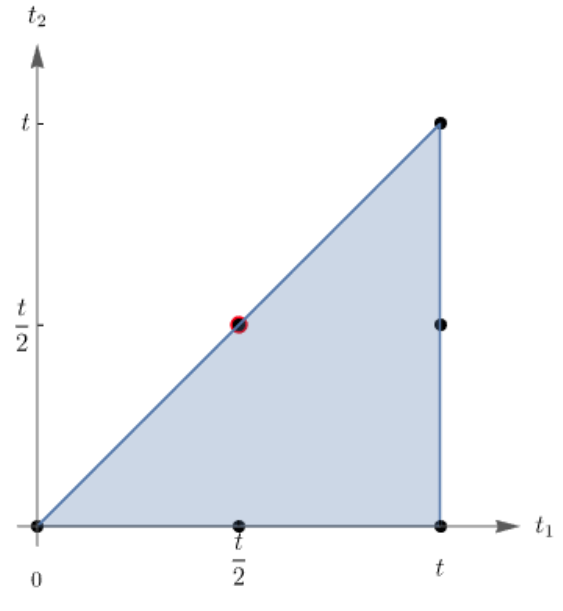}\hspace*{15pt}\includegraphics[width=200pt]{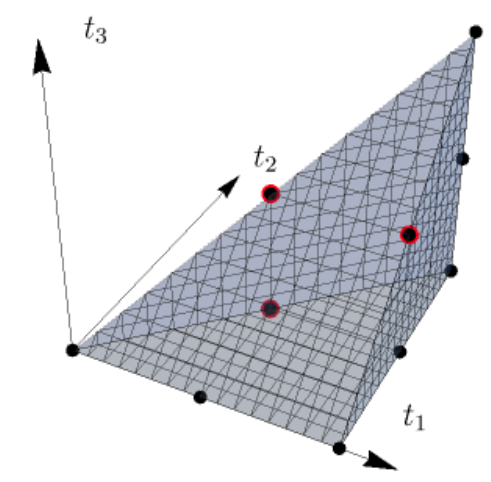}
    \caption{Two- and three-dimensional simplexes displaying the nodes for $\tau=0$.}
    \label{fig:2d-simplex}
  \end{center}
\end{figure}

Figure \ref{fig:2d-simplex} displays the relevant simplexes for $d=2,3$ with $\tau=0$ (the picture for general $\tau\in[0,\frac12)$ is similar and easy to envisage), highlighting the location of quadrature points. The function is evaluated at all the (small and large) discs, while directional derivatives are evaluated only at the large red discs. Note that {\em all\/} the quadrature points are of the form
\begin{displaymath}
  \mathbf{\alpha}_{d,\mathbf{j}}:=(\underbrace{\alpha_3,\ldots,\alpha_3}_{j_3 {\rm\ times}},\underbrace{\alpha_2,\ldots,\alpha_2}_{j_2 {\rm\ times}},\underbrace{\alpha_1,\ldots,\alpha_1}_{j_1{\rm\ times}} ), 
\end{displaymath}
where $j_1+j_2+j_3=d$ and $j_1,j_2,j_3\geq0$. Thus, for example, for $d=2$ we have six quadrature points
\begin{eqnarray*}
  &&(\alpha_1,\alpha_1):\quad \mathbf{j}=\left[
  \begin{array}{c}
  2\\0\\0
  \end{array}
  \right]\!,\qquad (\alpha_2,\alpha_1):\quad \mathbf{j}=\left[
  \begin{array}{c}
  1\\1\\0
  \end{array}
  \right]\!,\qquad (\alpha_3,\alpha_1):\quad \mathbf{j}=\left[
  \begin{array}{c}
  1\\0\\1
  \end{array}
  \right]\!,\\
  &&(\alpha_2,\alpha_2):\quad \mathbf{j}=\left[
  \begin{array}{c}
  0\\2\\0
  \end{array}
  \right]\!,\qquad (\alpha_3,\alpha_3):\quad \mathbf{j}=\left[
  \begin{array}{c}
  0\\0\\2
  \end{array}
  \right]\!,\qquad (\alpha_3,\alpha_2):\quad \mathbf{j}=\left[
  \begin{array}{c}
  0\\1\\1
  \end{array}
  \right]\!.
\end{eqnarray*}
Of them, just $\mathbf{j}^\top=\left[
  \begin{array}{ccc}
     0 & 2 & 0
  \end{array}
  \right]$ corresponds to a point denoted in Fig.~1 by a large disc, where we need to compute also derivatives.
  This  rule applies also to $d=3$: we need to compute derivatives only at points where $j_2\geq2$.\footnote{Of course, $j_2\geq2$ is impossible for $d=1$ where, indeed, we do not need to compute derivatives.}
  
  Denoting the $d$-dimensional simplex with vertices as in Fig.~1 by ${\mathcal S}_d$, 
\begin{equation*}
{\mathcal S}_d=\{(0,h)\times(0,\eta_1)\times...\times(0,\eta_{d-1})\subset \mathbf{R}^d\ |\ h\geq \eta_1\geq \ldots \geq \eta_{d-1}\geq0\},
\end{equation*} 
the quadrature rule used in this paper is thus
  \begin{equation}
    \label{quad}
    \int_{{\mathcal S}_d} f_d(\boldsymbol{\eta})\mathrm{d}\boldsymbol{\eta}= h^d \sum_{\mathbf{j}\in P_d} b_{d,\mathbf{j}} f_d(\mathbf{\alpha}_{d,\mathbf{j}})+h^{d+1} \sum_{\mathbf{j}\in Q_d} c_{d,\mathbf{j}} {\mathbf{v}_{d,\mathbf{j}}} ^\top \mathbf{\nabla}f_d(\mathbf{\alpha}_{d,\mathbf{j}})+R_{d},
  \end{equation}
  where $R_d=\mathcal{O}(h^5)$ as we will verify in Section \ref{section:ErrorTerms},
  \begin{displaymath}
    P_d=\{\mathbf{j}\in\mathbb{Z}_+^3\,:\, j_1+j_2+j_3=d\},\qquad Q_d=\{\mathbf{j}\in P_d,\:\, j_2\geq2\,\  \text{and}\ d=2,3\}
  \end{displaymath}
  and ${\mathbf{v}_{d,\mathbf{j}}}$ is  the direction of the inward normal at the point $\alpha_{d,\mathbf{j}}$.
  
  Note that $Q_1=Q_4=\emptyset$. 
  The differentiation directions are
      \begin{displaymath}
  \mathbf{v}_{2,{\tiny\left[
  \begin{array}{c}
  0\\2\\0
  \end{array}
  \right]\!}}=\left[
  \begin{array}{c}
         1\\-1
  \end{array}
  \right]\!,\quad\qquad \mathbf{v}_{3,{\tiny\left[
  \begin{array}{c}
  1\\2\\0
  \end{array}
  \right]\!}}=\left[
  \begin{array}{c}
         0\\1\\-1
  \end{array}
  \right]\!,\quad \mathbf{v}_{3,{\tiny\left[
  \begin{array}{c}
  0\\3\\0
  \end{array}
  \right]\!}}=\left[
  \begin{array}{c}
         1\\0\\-1
  \end{array}
  \right]\!,\quad \mathbf{v}_{3,{\tiny\left[
  \begin{array}{c}
  0\\2\\1
  \end{array}
  \right]\!}}=\left[
  \begin{array}{c}
         1\\-1\\0
  \end{array}
  \right]\!.
\end{displaymath}


The weights $b_{d,\mathbf{j}}$ and $c_{d,\mathbf{j}}$ are derived to ensure order 5. Specifically, we have observed that
\begin{equation}
  \label{weights}
  b_{d,\mathbf{j}}=\frac{v_1^{j_1}v_2^{j_2}v_3^{j_3}}{j_1!j_2!j_3!},
  \qquad 
  c_{2,\mathbf{j}}=
  \begin{cases}
      r\quad \mathbf{j}=[0\;2\;0]^\top\\
      0\quad {\rm otherwise},
  \end{cases} \quad
  c_{3,\mathbf{j}}=
  \begin{cases}
  v_3 r, \quad \mathbf{j}=[1\;2\;0]^\top\\
  \frac{1}{2!}v_2 r\quad \mathbf{j}=[0\;3\;0]^\top\\
  v_1 r\quad \mathbf{j}=[0\;2\;1]^\top\\
  0 \quad {\rm otherwise}.
\end{cases}
\end{equation}

While the quadrature (\ref{quad}) might seem at a first glance as unduly complicated -- after all, we require it just for $d=1,2,3,4$, while for $d=1$ and $d=4$ there are no derivative terms -- its form is conducive to further generalisation of our methodology to higher-order splittings. Explicit quadrature formulas are given in Appendix \ref{app:quadratures}.



This is the moment to remind the reader that, in general, we have $\tau\in[0,\frac12)$ and this has an obvious impact on both the location of cubature points and coefficients -- yet the principle is identical.

\section{Error Analysis}\label{section:ErrorTerms}
The error of the derived approximation (\ref{TACB}) originates in three sources. The first one comes from the truncation of the iterated Duhamel's formula. The second one follows from truncating the exponentials in (\ref{expansion}), noting that their arguments are bounded operators. The final one comes from the quadratures needed to integrate the iterated Duhamel's formula (\ref{iterations}).

\subsection{The first two error bounds}
This error bound for (\ref{RV}) can be obtained in a straightforward manner, bounding the integrand. It is convenient, however, to express the error bound of $R_{V}$ in terms of $\|u_0\|$, which is known. To do so, one can follow the approach  presented in \cite{FamilyStrang}. The latter results in the bound
\begin{equation}\label{eq:RV}
\|R_{V}\|\leq {\rm e}^{h\,C_h\,\|B\|_Y}\,\frac{(h\, C_h\,\|B\|_Y)^{5}}{5!}C_h\,\|u_0\|,
\end{equation}
  with the standard norm $\|v(\cdot)\|_Y:=\max_{s\in[0,h]}\|v(s)\|$ with $Y:=C^1([0,h],X)$. The constant $C_h$ depends on $h$ only, such that $\|{\rm e}^{tA}\|\leq C_h $ and $\|{\rm e}^{tB(\cdot)}\|_Y\leq C_h$  for $t\in[0,h]$. 
  
 The truncation of exponentials (\ref{expansion})   generates the error bounded by
  \begin{equation} \label{RE}
r_E=\frac{h^5}{120}(2p^5\|B\|^5+\|\tilde{B}\|^5)\|\mathrm{e}^{phA}\|^2\|\mathrm{e}^{(1-2p)hA}\|\|u_0\|
  \end{equation}
We recall that both $B$ and $\tilde{B}$ are bounded operators.
 
\subsection{Quadrature errors}

In our setting, the lion's share of the splitting error originates from the replacement of integrals by quadrature.
The order of accuracy of the quadratures can be derived by te standard approach: one simply has to determine the maximum total degree of  multivariate polynomials that are  integrated exactly by the quadrature rule. On the other hand, the determination of the error bounds requires further attention, especially in the case of multivariate integrals.

Quadrature errors will be written in terms of the (partial) derivatives of function $f_n$, see (\ref{fn}). An efficient way to compute them, allowing us to observe their structure, is provided by the the following lemma.
\begin{lemma}
\label{lemma:derivatives}
For sufficiently smooth $B(t)$ the $m$th derivative of $\ee^{-tA}B(t)\ee^{tA}$ with respect to $t$ can be expressed in the form
\begin{align*}
\partial_t^m(\ee^{-tA}B(t)\ee^{tA})&=\ee^{-tA}\sum_{k=0}^m\binom{m}{k}{\rm ad}_A^{(m-k)}B^{(k)}(t)\ee^{tA}\nonumber\\
&=\ee^{-tA}\sum_{k=0}^m\binom{m}{k}\underbrace{[[\ldots[[B^{(k)}(t),A],A]\ldots],A]}_{A \ \text{appears}\ m-k \ \text{times}}\ee^{tA}.
\end{align*}	

This implies that for $n=1,\ldots,4$ 
\begin{align*}
&\partial_{\eta_1}^{m_1}\partial_{\eta_2}^{m_2} \ldots \partial_{\eta_n}^{m_n}f_n(\eta_1,\eta_2,\ldots,\eta_n)
=f^{(m_1,m_2,\ldots,m_n)}_n(\eta_1,\eta_2,\ldots,\eta_n)=
\nonumber
\\
&\ee^{hA}\ee^{-\eta_1A}\sum_{k_1=0}^{m_1}\binom{m_1}{k_1}{\rm ad}_A^{m_1-k_1}B^{k_1}(t)\ee^{\eta_1A}\ee^{-\eta_2A}\sum_{k_2=0}^{m_2}\binom{m_2}{k_2}{\rm ad}_A^{m_2-k_2}B^{k_2}(t)\ee^{\eta_2A}
\cdots
\nonumber
\\
&\hspace*{15pt}\cdots\ee^{-\eta_nA}\sum_{k_n=0}^{m_n}\binom{m_n}{k_n}{\rm ad}_A^{m_n-k_n}B^{k_n}(t)\ee^{\eta_nA}u_0.
\end{align*}

\end{lemma}

This lemma is complemented with the following assumption that is required to derive and bound  quadratures errors.

\begin{assumption}\label{assum:derivatives}
	Let us assume that 
	\begin{align*}
		(a)\quad&f_1^{(4)},\\
		(b)\quad&f_2^{(i,j)},\quad i+j=3,\\
        (c)\quad&f_3^{(i,j,k)},\quad i+j+k=2,\\
        (d)\quad&f_4^{(i,j,k,l)},\quad i+j+k+l=1,
	\end{align*}
	are all well defined and bounded on $\mathcal{S}_d$, $d=1,\ldots,4$.
\end{assumption}

\subsubsection{The 1D Integral}

We start with the one-dimensional integral of (\ref{iterations}), which is calculated using the quadrature (\ref{1dimensional}). The error $R_1$, up to its dominant, term can be written explicitly using Taylor series with integral reminder term, see, e.g., \cite{powell81atm}. This  results in 
\begin{equation}
\label{eq:dominant}
R_1= \frac{ h^5 \left(10 \tau ^2-10 \tau +1\right)}{2880}f_1^{(4)}(0)+\mathcal{O}(h^6) . 
\end{equation}

It can be immediately verified that the dominant part of $R_1$ vanishes for $\tau=(5-\sqrt{15})/10$. This corresponds to Gauss--Legendre quadrature and  increases the order of the quadrature from $\mathcal{O}(h^5)$ to $\mathcal{O}(h^6)$-- needless to say,  only in the case of the 1D integral. 
 There is no other $\tau\in[0,1/2)$ with this property. This value of  $\tau$  defines the three-point Gauss--Legendre quadrature  in the interval $[0,h]$, whose nodes ($x$) satisfy $\mathrm{P}_3(2x/h-1)=0$, where $\mathrm{P}_3$ is the cubic Legendre polynomial.

In our particular case of interest, the fourth-order derivative required by (\ref{eq:dominant}) reads 
\begin{align}\label{eq:fourth_derivative}
f_1^{(4)}(0)=&\ee^{hA}\left\{[[[[B(0),A],A],A],A]+3[[[B'(0),A],A],A]+3[[B''(0),A],A]+[B'''(0),A]+B^{(4)}(0) \right\}u_0,
\end{align}
according to Lemma \ref{lemma:derivatives}.

\begin{figure}[H]
\begin{center}
\includegraphics[width=250pt]{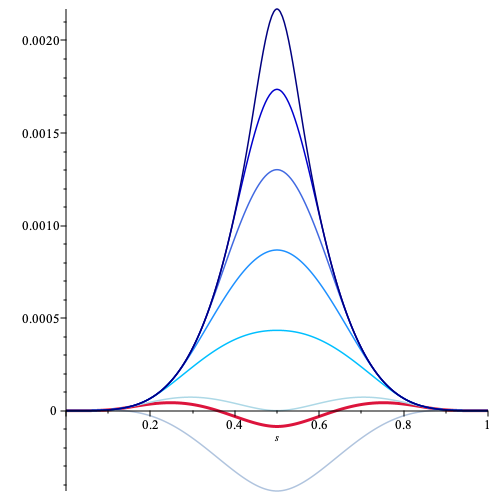}
\caption{The 1D Peano kernel (\ref{kernel}) for various values of $\tau$: $\tau=k/16$ for $k=1,\ldots,7$ in increasingly darker shades of blue and $\tau=(5-\sqrt{15})/10$ as a thick red line.}
\label{fig:2}
  \end{center}
\end{figure}

Alternatively, we can use the Peano kernel theorem \cite{powell81atm} to write $R_1$ in \textit{exact form} as
\begin{equation}
R_{1}\ =\ \int_0^{1}K(s;\tau)f_1^{(iv)}(hs)\,\text{d}s,
\end{equation}
where $K(s,\tau)$ is the {\em kernel\/} given by
\begin{equation}
\label{kernel}
K(s;\tau)\ =\ \frac{h^5}{24} \left\{
    \begin{aligned}
     s^4, &\qquad s \in[0,\tau],\\
      s^4-4p(s-\tau)^3, &\qquad s \in[\tau,\tfrac12],\\
     (1-s)^4 +4p(s-1+\tau)^3&\qquad s \in[\tfrac12,1-\tau],\\
      (1-s)^4, &\qquad s \in[1-\tau,1].
    \end{aligned}
    \right.
\end{equation}
Figure \ref{fig:2} shows the plot of the kernel for representative values of $\tau$. The explicit form of the kernel may be obtained using integration by parts, as done in \cite{FamilyStrang}.
From (\ref{kernel}), it is clear that the kernel $K(s,\tau)$ has a $B$-spline representation whose knot sequence is defined by the nodes.  Seeking $\tau$ that minimizes $|R_1|$ we recall that, by standard H\"older inequality,
\begin{equation}
	\label{eq:upperbounding_R1}
  |R_1|=\left|\int_0^1 K(s;\tau)f^{(iv)}(hs) \D s\right|\leq \|K(\,\cdot\,;\tau)\|_{\mathrm{L}_p[0,1]} \|f^{(iv)}\|_{\mathrm{L}_q[0,h]},
\end{equation}
where $p\in[1,\infty]$, $1/p+1/q=1$ \cite{powell81atm}. The two most pertinent values are $p=q=2$ and $p=1$, $q=\infty$ and in each case the optimal $\tau$ can be computed numerically:
\begin{equation}
	\label{eq:optimaltau_R1}
  \tau_2\approx 0.11760,\qquad \tau_\infty\approx 0.117882.
\end{equation}
Interestingly, both values are tantalisingly close to the unique zero of
\begin{displaymath}
\int_0^1 K(s;\tau)\,\text{d}s= \frac{h^5(10 \tau ^2-10 \tau +1)}{2880},
\end{displaymath}
in $(0,\frac12)$, namely
\begin{equation}
	\label{eq:trueoptimaltau_R1}
  \tau_{\mathrm{opt}}=\frac12-\frac{\sqrt{15}}{10}\approx 0.11270.
\end{equation}
Note that this value of $\tau$, which increases the order to ${\cal O}(h^6)$, allows for an enhanced Peano kernel estimate, but this is of little significance in our setting because of higher-dimensional integral.

Yet, it is clear from the above discussion that $|R_1|$ displays very little sensitivity to the values of $\tau$ in the vicinity of $\tau_{\mathrm{opt}}$ and, for simplicity's sake, it is a good idea to use the latter as a working approximation to the optimal value of the parameter -- at least, for the time being,  in the univariate case.


\subsubsection {The 2D Integral}

The error of the two-dimensional integral $R_2$ can be analysed in the same manner as for the one-dimensional one. First, we establish the dominant term of $R_2$ using Taylor series, namely 
\begin{align}\label{eq:R2}
R_2=& -\frac{h^5}{8640 (2 \tau -1)} 
\left\{
\left(
120 \tau ^2-84 \tau +7\right) [f_2^{(0,3)}(0,0)-f_2^{(3,0)}(0,0)]
\right.
\nonumber \\
&
\left.
+(36 \tau -3) [f_2^{(1,2)}(0,0)- f_2^{(2,1)}(0,0)]
\right\}
+\mathcal{O}(h^6). 
\end{align}
where
\begin{equation*}
f_2^{(0,3)}(0,0)-f_2^{(3,0)}(0,0)=\ee^{hA}[B(0),[[[B(0),A],A],A]+3[[B'(0),A],A]+3[B''(0),A]+B'''(0)]u_0,
\end{equation*}
\begin{equation*}
f_2^{(1,2)}(0,0)- f_2^{(2,1)}(0,0)=\ee^{hA}[[B(0),A]+B'(0),[[B(0),A],A]+2[B'(0),A]+B''(0)]u_0,
\end{equation*}
according to Lemma  \ref{lemma:derivatives}, because
\begin{align}
f_2^{(3,0)}(0,0)=&\ee^{hA}\left\{[[[B(0),A],A],A]+3[[B'(0),A],A]+3[B''(0),A]+B'''(0) \right\}B(0),
\end{align}
\begin{align}
f_2^{(0,3)}(0,0)=&\ee^{hA}B(0)\left\{[[[B(0),A],A],A]+3[[B'(0),A],A]+3[B''(0),A]+B'''(0) \right\}
\end{align}
\begin{align}
f_2^{(2,1)}(0,0)=&\ee^{hA}\left\{[[B(0),A],A]+2[B'(0),A]+B''(0) \right\}\left\{[B(0),A]+B'(0)\right\},
\end{align}
\begin{align}
f_2^{(1,2)}(0,0)=&\ee^{hA}\left\{[B(0),A]+B'(0)\right\}\left\{[[B(0),A],A]+2[B'(0),A]+B''(0) \right\},
\end{align}

Consider the expressions
\begin{displaymath}
   P_a(\tau):=-\frac{120\tau^2-84\tau+7}{8640(2\tau-1)}\qquad\mbox{and}\qquad P_b(\tau):=\frac{36\tau-3}{8640(2\tau-1)},
\end{displaymath}
appearing in (\ref{eq:R2}). Their values at $\tau_{\rm{opt}}\approx 0.11270$ are
\begin{displaymath}
   P_a(\tau_{\rm opt})\approx-0.00014085\qquad\mbox{and}\qquad P_b(\tau_{\rm opt})\approx-0.000157967,
\end{displaymath}

From (\ref{eq:R2}) it can be checked that there is no value of $\tau\in[0,\frac12)$ that increases the order of $R_2$ to $O(h^6)$. 

As was the case for 1D, further insight may be gleaned if we look at the exact form of the remainder $R_2$.
To do so,  we use the concept of Sard kernels, see  \cite{Engels1980,sard63la}, an extension of Peano kernel theory to a multivariate setting. In this framework and following much algebra, the error is written as
\begin{eqnarray}
\label{eq:sard}	
R_2&\!\!\!=\!\!\!&
\int_0^1 \{K_{3,0}(s;\tau) f_2^{(3,0)}(s,0) + K_{0,3}(s;\tau) f_2^{(0,3)}(0,s)\}\D s \\ \nonumber
&\!\!\!\!\!\!&\mbox{}
+ \!\int_0^1 \{K_{2,1}(s;\tau) f_2^{(2,1)}(s,0)+ K_{1,2}(s;\tau) f_2^{(1,2)}(0,s)\}\D s\,\\ \nonumber
&\!\!\!\!\!\!&\mbox{}
+\int_0^1 \!\int_0^1 \{K_{2,1}(s_1,s_2;\tau)f_2^{(2,1)}(s_1,s_2)+
K_{1,2}(s_1,s_2;\tau) f_2^{(1,2)}(s_1,s_2)\}\D s_2\D s_1\\ \nonumber
&\!\!\!\!\!\!&\mbox{}
+\int_0^1\int_0^1  K_{2,1}^{\delta}(s_1,s_2;\tau)f_2^{(2,1)}(s_1,s_2)\D s_1\D s_2
	+
	\int_0^1\int_0^1  K_{1,2}^\delta(s_1,s_2;\tau)f_2^{(1,2)}(s_1,s_2)\D s_2
\end{eqnarray}
For the quadratures used in this study, related Sard kernels present a cumbersome structure, especially those which are $(s_1,s_2)$-dependent. However, they can be easily handled numerically. That is why, for the sake of simplicity of the presentation, we present  kernels $K_{3,0}$, $K_{0,3}$, $K_{2,1}^\delta(s_1,s_2;\tau)$, and $K_{2,1}^\delta(s_1,s_2;\tau)$   only. Their expressions are 

\begin{gather*}
K_{3,0}(s;\tau)	
=\frac{h^5}{24} (s-1)^3 (s+3)
\\
+\frac{h^5}{2} \begin{cases}
	0,& 1-\tau\leq s \leq 1 \\
		\left(\dfrac{v_3^2}{2!}+v_3v_2+v_3v_1\right)(s+\tau-1)^2,&\tfrac{1}{2}\leq s\leq 1-\tau,\\[4pt]
		\left(\dfrac{v_3^2}{2!}+v_3v_2+v_3v_1\right)(s+\tau-1)^2+\left(\dfrac{v_2^2}{2!}+v_2v_1\right)(s-\tfrac{1}{2})^2-r(2s-1),&\tau\leq s\leq\tfrac{1}{2},\\[4pt]
		\left(\dfrac{v_3^2}{2!}+v_3v_2+v_3v_1\right)(s+\tau-1)^2+\left(\dfrac{v_2^2}{2!}+v_2v_1\right)(s-\tfrac{1}{2})^2+\dfrac{v_1^2 }{2!}(s-\tau)^{2},\\
		\hspace*{20pt}\mbox{}-r(2s-1),& 0\leq s\leq \tau;
	\end{cases}
\end{gather*}	

\begin{gather*}
	K_{0,3}(s;\tau)	
	=-\frac{h^5}{24} (s-1)^4
	\\
	+\frac{h^5}{2}\begin{cases}
		0,& 1-\tau\leq s\leq 1, \\[4pt]
		\dfrac{v_3^2}{2!}(s+\tau-1)^2,&\frac{1}{2}\leq s\leq1-\tau,\\[4pt]
		\dfrac{v_3^2}{2!}(s+\tau-1)^2+\left(v_3v_2+\dfrac{v_2^2}{2!}\right)(s-\frac{1}{2})^{2}+r(2s-1),&\tau\leq s\leq \frac{1}{2},\\[4pt]
		\dfrac{v_3^2}{2!}(s+\tau-1)^2+\left(v_3v_2+\dfrac{v_2^2}{2!}\right)(s-\frac{1}{2})^{2}+\left(v_3v_1+v_2v_1+\dfrac{v_1^2 }{2!}\right)(s-\tau)^{2}\\
		\hspace*{20pt}\mbox{}+r(2s-1),& 0 \leq s\leq\tau;
	\end{cases}
\end{gather*}	
and
\begin{align}
K_{2,1}^\delta(s_1,s_2;\tau)&=h^5r(2s_1-1)\delta(s_2-1)(\tfrac{1}{2}-s_1)_+ ,\nonumber\\ K_{1,2}^\delta(s_1,s_2;\tau)&=h^5r(2s_2-1)\delta(s_1-1)(\tfrac{1}{2}-s_2)_+,
\end{align}	
where $(x)_+=\max\{x,0\}$. 
The kernels labeled by $\delta$ stand out due to the present of a Dirac delta. All other kernels do not contain this distribution. As mentioned before, the other kernels can be handled numerically. Since there is no value that increases the order of $R_2$, one may find the optimal value for $\tau$ minimizing $|R_2|$. Due to its structure, see (\ref{eq:sard}), we minimise the  each Sard kernel as in (\ref{eq:upperbounding_R1}). In Table \ref{table:optimal_R2}, we show the values for optimal $\tau$ in $\mathrm{L}_1$ and $\mathrm{L}_2$ norms.

\begin{table}[h]
	\centering
	\caption{Values for $\tau$ that minimise the Sard kernels in the  $\mathrm{L}_1$ and $\mathrm{L}_2$ norms, see text.}
		{\setlength{\tabcolsep}{0.4cm}	
	\begin{tabular}{ccc}
		\hline
		\rule{-2pt}{4ex}Kernel              & $\ta_\infty$ & $\tau_2$  \\[5pt]
		\hline
		\rule{-2pt}{4ex}
		$K_{3,0}(s;\tau)$, $K_{0,3}(s;\tau)$ &    0.110387          &   0.113077                     \\[5pt]
		$K_{2,1}(s;\tau)$, $K_{1,2}(s;\tau)$& 0.119818             & 0.11757                              \\[5pt]
		$K_{2,1}(s_1,s_2;\tau)$, $K_{1,2}(s_1,s_2;\tau)$&    0.1          &  0.1          \\[5pt]
		$K_{2,1}^\delta(s_1,s_2;\tau)$,  $K_{1,2}^\delta(s_1,s_2;\tau)$&  0.146447            &     0.146447      \\[10pt]
		\hline                    
	\end{tabular}
}
	\label{table:optimal_R2}
\end{table}
Except for the kernels $K_{2,1}^\delta$ and $K_{1,2}^\delta $, the values shown in Table \ref{table:optimal_R2} are close to those for the univariate quadratures shown in (\ref{eq:optimaltau_R1}) and also close to the enhanced value $\tau_{\text{opt}}$ shown in (\ref{eq:trueoptimaltau_R1}).


\subsubsection{The 3D and 4D Integrals}
In our error analysis of  three- and four-dimensional quadratures, we restrict ourselves to  Taylor expansions with integral remainder. While in principle it is possible to compute  Sard kernels for these quadratures, their expressions are exceedingly cumbersome. For the particular purpose of the work, the connection between the kernels and the dominant error terms provides enough sufficient insight  to forego the analysis of Sard kernels.
For the three-dimensional quadrature, we have
\begin{align}\nonumber
R_3=& -\frac{h^5 \left(432 \tau ^4-864 \tau ^3+558 \tau ^2-126 \tau +7\right)}{12960 (1-2 \tau )^4}[f_3^{(2,0,0)}(0,0,0)-2f_3^{(0,2,0)}(0,0,0)+f_3^{(0,0,2)}(0,0,0)] \\ \nonumber
&\frac{h^5 \left(216 \tau ^4-252 \tau ^3+99 \tau ^2-18 \tau +1\right)}{3240 (1-2 \tau )^4}[f_3^{(1,1,0)}(0,0,0)-2f_3^{(1,0,1)}(0,0,0)+f_3^{(0,1,1)}(0,0,0)]
\\ \label{eq:R3}
&+\mathcal{O}(h^6) 
\end{align}    
where, according to Lemma \ref{lemma:derivatives}, 
$$
f_3^{(2,0,0)}(0,0,0)-2f_3^{(0,2,0)}(0,0,0)+f_3^{(0,0,2)}(0,0,0)=\ee^{hA}[[[[B(0),A],A]+2[B'(0),A]+B''(0),B(0)],B(0)]u_0
$$
and
$$
f_3^{(1,1,0)}(0,0,0)-2f_3^{(1,0,1)}(0,0,0)+f_3^{(0,1,1)}(0,0,0)=\ee^{hA}[[B(0),A]+B'(0),[[B(0),A]+B'(0),B(0)]]u_0.
$$

Let us have a more careful look at the size of 
\begin{displaymath}
   Q_a(\tau):=-\frac{432 \tau ^4-864 \tau ^3+558 \tau ^2-126 \tau +7}{12960 (1-2 \tau )^4}\qquad\mbox{and}\qquad Q_b(\tau):=\frac{216 \tau ^4-252 \tau ^3+99 \tau ^2-18 \tau +1}{3240 (1-2 \tau )^4},
\end{displaymath}
the coefficients in (\ref{eq:R3}). For $\tau\in[0,\frac12]$ the first vanishes at $1/2-\sqrt{15\pm\sqrt{105}}/12$, the latter at $\approx 0.0922146$. Recalling, however, that we have recommended in Subsection~3.2.1 to use $\tau_{\mathrm{opt}}=1/2-\sqrt{15}/10$, we note that
\begin{displaymath}
  Q_a(\tau_{\mathrm{opt}})=\frac{1}{33645}\approx 0.0002743,\qquad Q_b(\tau_{\mathrm{opt}})=\frac{259}{29160}-\frac{\sqrt{14}}{432}\approx -0.0000832,
\end{displaymath}
both fairly small.

In turn, the error of the integral over the four-dimensional simplex reads
\begin{align}\nonumber
R_4=&\  h^5R_a(\tau)[f_4^{(1,0,0,0)}(0,0,0,0)-f_4^{(0,0,0,1)}(0,0,0,0)]+h^5R_b(\tau)[f_4^{(0,1,0,0)}(0,0,0,0)-f_4^{(0,0,1,0)}(0,0,0,0)]\\ \label{eq:R4}
&+\mathcal{O}(h^6) 
\end{align} 
where
\begin{align*}
R_a(\tau)&=\frac{248832 \tau ^7-732672 \tau ^6+891648 \tau ^5-578880 \tau ^4+216000 \tau ^3-46416 \tau ^2+5376 \tau -269}{155520 (2 \tau -1)^7} \\
R_b(\tau)&=\frac{27648 \tau ^7-96768 \tau ^6+145152 \tau ^5-118080 \tau ^4+54720 \tau ^3-14064 \tau ^2+1824 \tau -91}{51840 (2 \tau -1)^7}
\end{align*}
and
$$
f_4^{(1,0,0,0)}(0,0,0,0)-f_4^{(0,0,0,1)}(0,0,0,0)=e^{hA}[[B(0),A]+B'(0),B^3(0)]u_0
$$
$$
f_4^{(0,1,0,0)}(0,0,0,0)-f_4^{(0,0,1,0)}(0,0,0,0)=e^{hA}[B(0)\{[B(0),A]+B'(0)\}B(0),B(0)]u_0
$$
There is no value of $\tau$ that increases the order of both $R_3$ and $R_4$ simultaneously to $O(h^6)$. Moreover,
\begin{displaymath}
  R_a(\tau_{\mathrm{opt}})=\frac{1}{80}-\frac{7579\sqrt{15}}{2519424}\approx 0.0008492,\qquad R_b(\tau_{\mathrm{opt}})=\frac{1}{240}-\frac{845\sqrt{15}}{839808}\approx 0.0002697,
\end{displaymath}
both fairly small.


 Based on our analysis, we  state the following Theorem on the convergence of family $\mathcal{T}_{ACB}^{(4)}$:
\begin{thm}
	Under Assumptions \ref{assum:main_assum}, \ref{assum:commutators} and  \ref{assum:derivatives}, the family of integrators (\ref{TACB})  satisfies fourth order convergence with

		\begin{equation} 
		\|u(h)- \mathcal{T}_{ACB}^{(4)}u_0\| \leq \|R_V\|+r_E+\sum_{i=1}^4\|R_i\|,  \label{estimate1}
	\end{equation}
	where constants  $\|R_V\|$, $r_E$ and $\|R_i\|$, $i=1,\ldots,4$ have been given by (\ref{eq:RV}), (\ref{RE}) and (\ref{eq:dominant}), (\ref{eq:R2}), (\ref{eq:R3}) and (\ref{eq:R4}), respectively.
\end{thm}

The proof of the Theorem follows directly from the estimates of the error terms derived in the previous Subsections. 

We reiterate that no single value of $\tau\in[0,\frac12]$ minimises all the remainders $R_i$, $i=1,2,3,4$ simultaneously. Assuming that the most nested commutator brings the biggest error, we wish to eliminate term (\ref{eq:fourth_derivative}) from the error $R_1$ defined with (\ref{eq:dominant}). This can be done setting $\tau=\tau_{\mathrm{opt}}$. As we have seen, this choice leads to small error constants in other integrals.

\section{A numerical example}\label{section:numerical_example}
Based on the  one-dimensional linear Schr\"odinger equation, we illustrate the analysis of this paper by numerical experiments. We consider the following  equation 
\begin{equation*}\label{Schrodinger}
	\left\{
	\begin{aligned}
		&\mathrm{i}\frac{\partial u(x,t)}{\partial t} =\left(-\frac{1}{2}\frac{\partial^2 u(x,t)}{\partial x^2} + V(x,t)\right)u(x,t),\qquad t\in[0,1],\qquad x\in[-\infty,\infty],\\
		&u(x,0)\ =\ \frac{\sin (20 (x-3))}{1+x^{10}},\\
		&u(\pm \infty,t)=0,
	\end{aligned}
	\right.
\end{equation*}
with the potential
\begin{equation*}
	V(x,t) =  (x-t)^2.
\end{equation*}

A finite-difference scheme of order \(\mathcal{O}(\Delta x^8)\) with \(\Delta x = 0.008\) ensures high spatial accuracy in the numerical solution. Then, using time-stepping on \(\mathcal{T}_{ACB}^{(4)}\), we evolved the initial condition from \(t = 0\) to \(t = 1\). The domain has been restricted to \(x \in [-40, 40]\), sufficiently large to avoid reflections due to  artificial boundaries. At the end of the evolution, we compared the numerical solution with the reference solution (obtained with a time step \(h = 10^{-4}\)) using the \(L^2[-\infty, \infty]\) norm. 

\begin{figure}[h]
	\centering
	\includegraphics[scale=0.2]{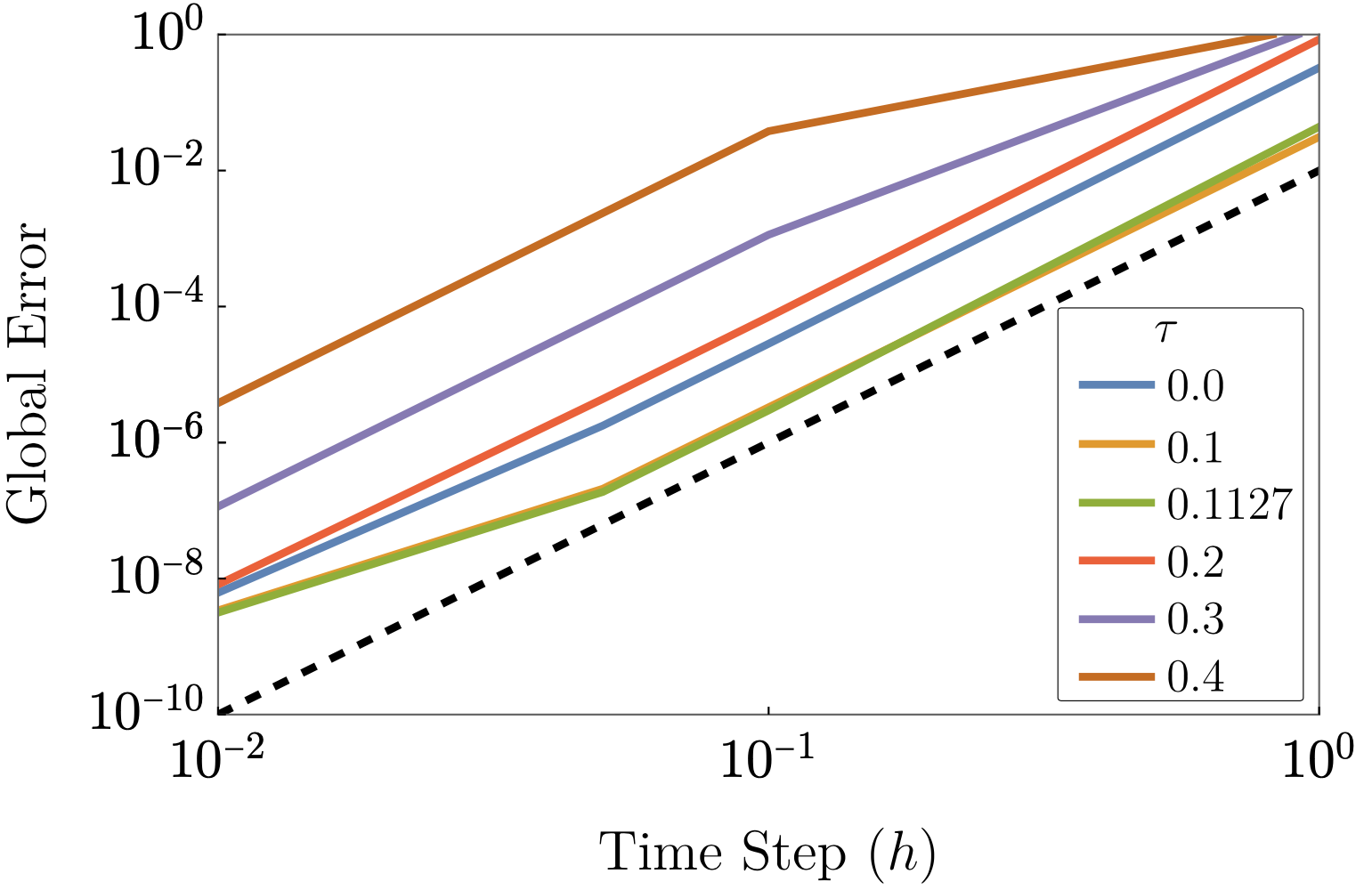}
	\caption{Global Error (in log-log scale) for representative values of $\tau$ as a function of the time step $h$. The black-dashed line represents the plot of $h^4$. The values $\tau=0.1127$ and $0.1$ results in similar global errors. }
	\label{fig:error}
\end{figure}

Figure \ref{fig:error} shows the global error as a function of \(h\) for representative values of \(\tau \in [0, 1/2)\).

We observe that the global error increases as one considers values of $\tau$ closer to $\tau=\frac12$. This is expected since  $R_i$, $i=2,3,4$, become large when $\tau\rightarrow\frac12$.

Consistently with our theoretical results related to the optimal member of the family, we observe that $\tau=\tau_{\mathrm{opt}} \approx \frac{1}{10}$ leads to a smaller error compared to the other values considered.  
Based on a numerical example, the authors of \cite{Chin} (the work where the family $\mathcal{T}_{ACB}^{(4)}$ was proposed) looked for the value of optimal $\tau$ scanning  values $\tau=0.142,0.144,0.146$. In their exploration, they found that \(\tau \approx 0.144\) leads to the smallest error. To some extent, this value corroborate our findings, namely $\tau_{opt}=0.11270$.


\section{Conclusions}\label{section:conclusions}

Compact splittings, as introduced in \cite{Chin}, represent a powerful means to solve the linear Schr\"odinger equation with variable potential in a precise and efficient manner. Yet, once the underlying domain, hence also the Laplace operator, are unbounded, a rigorous proof of their convergence is missing. In this paper we have presented such a proof. Moreover, the splittings include a parameter, whose choice has significant impact on the size of the error. Providing careful error analysis, we argue that, while a single choice of the parameter $\tau$ cannot minimise simultaneously all the error components, an excellent choice is $\tau_{\mathrm{opt}}=1/2-\sqrt{15}/10$. This is affirmed by a numerical experiment. A good choice of the parameter is important: Fig.~\ref{fig:error} demonstrates that it can lead to significant improvement in the quality of the solution, up to orders of magnitude,  while not imposing any extra cost.

\section{Acknowledgements}
The work of Karolina Kropielnicka and Juan Carlos Del Valle Rosales in this project was funded by the National Science Centre, Poland, (NCN) project no. 2019/34/E/ ST1/00390.

\bibliographystyle{abbrv}
\bibliography{main1}

\appendix

\def\theequation{\Alph{section}.\arabic{equation}}
\section{Explicit quadrature formul\ae}
\label{app:quadratures}
We recall from (\ref{1dimensional}) that
\begin{displaymath}
  \int_{\mathcal{S}_1} f_1(\eta_1)\text{d}\eta_1 =h \sum_{k=1}^3 v_{k} f_1(\alpha_k) + R_{1}.
\end{displaymath}
Likewise -- here and in the sequel the weights $b_{d,\mathbf{j}}$ are consistent with (\ref{weights}) --
\begin{eqnarray}\label{2dimensional}
        &\!\!\!\!\!\!&\int_{\mathcal{S}_2} f_2(\eta_1,\eta_2)\,\text{d}\eta_2\text{d}\eta_1\\
        \nonumber
        &\!\!\!=\!\!\!& h^2\sum_{\mathbf{j}\in P_2} b_{2,\mathbf{j}}f_2(\mathbf{\alpha}_{2,\mathbf{j}})+h^3 r \mathbf{v}_{1,{\tiny\left[
  \begin{array}{c}
  0\\2\\0
  \end{array}
  \right]}}^\top\grad f_2(\alpha_2,\alpha_2)+R_2\\
        \nonumber
        &\!\!\!=\!\!\!&h^2\left[\frac{v_3^2}{2!}f_2(\alpha_3,\alpha_3)+v_3v_2f_2(\alpha_3,\alpha_2)+v_3v_1f_2(\alpha_3,\alpha_1)+\frac{v_2^2}{2!}f_2(\alpha_2,\alpha_2)+v_2v_1f_2(\alpha_2,\alpha_1)+\frac{v_1^2}{2!}f_2(\alpha_1,\alpha_1)\right]\\
        \nonumber
        &\!\!\!\!\!\!&\mbox{}+h^3r\, [\partial_{\eta_1}f_2(\alpha_2,\alpha_2) -\partial_{\eta_2}f_2(\alpha_2,\alpha_2)]+R_2,
\end{eqnarray} 
where $R_2=\mathcal{O}(h^5)$, and 
\begin{eqnarray}
\label{3dimensional}
  &\!\!\!\!\!\!&\int_{\mathcal{S}_3} f_3(\eta_1,\eta_2,\eta_3)\,\text{d}\eta_3\text{d}\eta_2\text{d}\eta_1=h^3 \sum_{\mathbf{j}\in P_3} b_{3,\mathbf{j}}f_3(\alpha_{3,\mathbf{j}}) \\ \nonumber
 &\!\!\!\!\!\!&\mbox{}+h^4 v_3\,r\, \mathbf{v}_{3,{\tiny\left[
  \begin{array}{c}
  1\\2\\0
  \end{array}
  \right]}}^\top \grad f_3(\alpha_3,\alpha_2,\alpha_2) +h^4 \frac{v_2\,r}{2!}\mathbf{v}_{3,{\tiny\left[
  \begin{array}{c}
  0\\3\\0
  \end{array}
  \right]}}^\top \grad f_3(\alpha_2,\alpha_2,\alpha_2)+h^4 v_1\,r\, \mathbf{v}_{3,{\tiny\left[
  \begin{array}{c}
  0\\2\\1
  \end{array}
  \right]}}^\top \grad f_3(\alpha_2,\alpha_2,\alpha_1)+R_3 \nonumber \\ 
 &\!\!\!=\!\!\!&h^3 \left[\frac{v_3^3}{3!}f_3(\alpha_3,\alpha_3,\alpha_3)+\frac{v_3^2v_2}{2!} f_3(\alpha_3,\alpha_3,\alpha_2)+\frac{v_3^2v_1}{2!}f_3(\alpha_3,\alpha_3,\alpha_1)+\frac{v_3v_2^2}{2!}f_3(\alpha_3,\alpha_2,\alpha_2)+v_3v_2v_1f_3(\alpha_3,\alpha_2,\alpha_1) \right.\nonumber \\ 
 &\!\!\!\!\!\!&\left.\hspace*{20pt}\mbox{}+\frac{v_2^3}{3!}f_3(\alpha_2,\alpha_2,\alpha_2)+\frac{v_2^2v_1}{2!}f_3(\alpha_2,\alpha_2,\alpha_1)+\frac{v_3v_1^2}{2!} f_3(\alpha_3,\alpha_1,\alpha_1)+\frac{v_3v_1^2}{2!} f_3(\alpha_3,\alpha_1,\alpha_1)+\frac{v_1^3}{3!}f_3(\alpha_1,\alpha_1,\alpha_1)\right] \nonumber \\
 \nonumber
 &\!\!\!\!\!\!&\mbox{}+h^4 pr [\partial_{\eta_2} f_3(\alpha_3,\alpha_2,\alpha_2)-\partial_{\eta_3} f_3(\alpha_3,\alpha_2,\alpha_2)]+\frac{h^4q\,r}{2!} [\partial_{\eta_1} f_3(\alpha_2,\alpha_2,\alpha_2)-\partial_{\eta_3} f_3(\alpha_2,\alpha_2,\alpha_2)]\\
 &\!\!\!\!\!\!&\mbox{}+h^4pr [\partial_{\eta_1} f_3(\alpha_2,\alpha_2,\alpha_1)-\partial_{\eta_2} f_3(\alpha_2,\alpha_2,\alpha_1)]+R_3
 \end{eqnarray}	
where $R_3=\mathcal{O}(h^5)$.
Finally,
\begin{eqnarray}
  \label{4dimensional}
  &\!\!\!\!\!\!&\int_{\mathcal{S}_4} f_4(\eta_1,\eta_2,\eta_3,\eta_4)\,\text{d}\eta_4\text{d}\eta_3\text{d}\eta_2\text{d}\eta_1\\
  \nonumber
  &\!\!\!=\!\!\!&h^4\sum_{\mathbf{j}\in P_4} b_{4,\mathbf{j}} f_4(\alpha_{4,\mathbf{j}})
 + \mathcal{O}(h^5)\\
 \nonumber
 &\!\!\!=\!\!\!&h^4[\frac{v_3^4}{4!}f_4(\alpha_3,\alpha_3,\alpha_3,\alpha_3)+\frac{v_3^3v_2}{3!} f_4(\alpha_3,\alpha_3,\alpha_3,\alpha_2)+\frac{v_3^3v_1}{3!}f_4(\alpha_3,\alpha_3,\alpha_3,\alpha_1)+\frac{v_3^2v_2^2}{2!2!}f_4(\alpha_3,\alpha_3,\alpha_2,\alpha_2)\\
 \nonumber
 &\!\!\!\!\!\!&\hspace*{20pt}\mbox{}+
 \frac{v_3^2v_2v_1}{2!} f_4(\alpha_3,\alpha_3,\alpha_2,\alpha_1) +\frac{v_3^2v_1^2}{2!2!}f_4(\alpha_3,\alpha_3,\alpha_1,\alpha_1)+\frac{v_3v_2^3}{3!}f_4(\alpha_3,\alpha_2,\alpha_2,\alpha_2)+\frac{v_3v_2^2v_1}{2!} f_4(\alpha_3,\alpha_2,\alpha_2,\alpha_1)\\
 \nonumber
 &\!\!\!\!\!\!&\hspace*{20pt}\mbox{}+\frac{v_3v_2v_1^2}{2!}f_4(\alpha_3,\alpha_2,\alpha_1,\alpha_1)+\frac{v_3v_1^3}{3!}f_4(\alpha_3,\alpha_1,\alpha_1,\alpha_1)+\frac{v_2^4}{4!} f_4(\alpha_2,\alpha_2,\alpha_2,\alpha_2)+\frac{v_2^3v_1}{3!} f_4(\alpha_2,\alpha_2,\alpha_2,\alpha_1)\\
 \nonumber
 &\!\!\!\!\!\!&\hspace*{20pt}\mbox{}+\frac{v_2^2v_1^2}{2!2!}f_4(\alpha_2,\alpha_2,\alpha_1,\alpha_1)+\frac{v_2^2v_1^3}{2!3!} f_4(\alpha_2,\alpha_1,\alpha_1,\alpha_1)+\frac{v_1^4}{4!}f_4(\alpha_1,\alpha_1,\alpha_1,\alpha_1)]+R_4
\end{eqnarray} 
and $R_4=\mathcal{O}(h^5)$.

\section{Explicit terms of the type (\ref{NiceTerm})}\label{explicitterms}

We present here a list of all 35 products of terms that enter into the expression (\ref{NiceTerm}).
\begin{enumerate}
\item[]  {\bf $\mathcal{O}(1)$ terms}
\item
$\ee^{hA}u_0$

\item[]  {\bf $\mathcal{O}(h)$ terms}
\item
$hp\, \ee^{h\tau A} \times B((1-\tau)h)  \times  \ee^{h(1-\tau)A}u_0\ 
\leadsto hv_3f_1(\alpha_3)$
 \item
$hq\,\ee^{\frac{1}{2} h A} \times  B(\tfrac12 h)  \times \ee^{\frac{1}{2} h A}u_0\ 
\leadsto  hv_2f_1(\alpha_2)$
 \item
 $hp\,\ee^{h(1-\tau)A} \times B(\tau h) \times \ee^{h\tau A}u_0\ 
 \leadsto  hv_1f_1(\alpha_1)$

\item[] {\bf $\mathcal{O}(h^2)$ terms}
\item
$\frac12 h^2 p^2 \ee^{h\tau A}  \times  B^2((1-\tau)h)  \times  \ee^{h(1-\tau)A}u_0\ 
\leadsto \frac{1}{2}h^2v_3^2f_2(\alpha_3,\alpha_3)$
 \item
$\frac12h^2 q^2 \ee^{\frac{1}{2} h A} \times   B^2(\tfrac12 h)  \times \ee^{\frac{1}{2} h A}u_0\ 
\leadsto \frac{1}{2}h^2v_2^2f_2(\alpha_2,\alpha_2)$%
  \item
 $\frac12 h^2 p^2 \ee^{h(1-\tau)A} \times B^2(\tau h) \times \ee^{h\tau A}u_0\ 
 \leadsto \frac{1}{2}h^2v_1^2f_2(\alpha_1,\alpha_1)$%
 \item
 $h^2 pq\,  \ee^{h\tau A} \times  B((1-\tau)h) \times \ee^{h(\frac12-\tau)A} \times B(\tfrac12 h) \times \ee^{\frac12hA}\ \leadsto h^2v_3v_2f_2(\alpha_3,\alpha_2)$
 \item
$h^2 pq\, \ee^{h\frac12 A} \times B(\tfrac12 h) \times \ee^{h(\frac12-\tau)A} \times B(\tau h) \times \ee^{h\tau A}u_0\   \leadsto h^2v_2v_1f_2(\alpha_2,\alpha_1) $%
\item  
$h^2 p^2 \ee^{h\tau A} \times B((1-\tau)h)  \times \ee^{h(1-2\tau)A} \times B(\tau h) \times \ee^{h\tau A}u_0\ 
\leadsto h^2v_3v_1f_2(\alpha_3,\alpha_1) $

\item[] {\bf $\mathcal{O}(h^3)$ terms}

 \item
$\frac16 h^3 p^3 \ee^{h\tau A} \times  B^3((1-\tau)h)  \times  \ee^{h(1-\tau)A}u_0\ \leadsto \frac{1}{3!}h^3v_3^3f_3(\alpha_3,\alpha_3,\alpha_3)$%
  \item
 $\frac16 h^3 p^3 \ee^{h(1-\tau)A} \times B^3(\tau h) \times \ee^{h\tau A}u_0\ 
 \leadsto \frac{1}{3!}h^3v_1^3f_3(\alpha_1,\alpha_1,\alpha_1)$%
  \item
 $ \frac12 h^3 pq^2 \ee^{h\tau A} \times  B((1-\tau)h) \times \ee^{h(\frac12-\tau)A} \times  B^2(\tfrac12 h) \times \ee^{\frac12hA}u_0\ 
  \leadsto \frac{1}{2!}h^3v_3v_2^2f_3(\alpha_3,\alpha_2,\alpha_2)$%
\item
$ \frac12 h^3 pq^2\ee^{h\frac12 A} \times B^2(\tfrac12 h) \times \ee^{h(\frac12-\tau)A} \times B(\tau h) \times \ee^{h\tau A}u_0\ 
  \leadsto \frac{1}{2!}h^3v_2^2v_1f_3(\alpha_2,\alpha_2,\alpha_1)$%
 \item
$\frac12 h^3 p^2q\, \ee^{h\frac12 A} \times B(\tfrac12 h) \times \ee^{h(\frac12-\tau)A} \times  B^2(\tau h) \times \ee^{h\tau A}u_0\ 
\leadsto \frac{1}{2!}h^3v_2v_1^2f_3(\alpha_2,\alpha_1,\alpha_1)$%
  \item
  $\frac12 h^3 p^2 q\,\ee^{h\tau A} \times  B^2((1-\tau)h) \times  \ee^{h(\frac12-\tau)A} \times B(\tfrac12 h) \times \ee^{h\frac12A}u_0\ 
  \leadsto \frac{1}{2!}h^3v_3^2v_2f_3(\alpha_3,\alpha_3,\alpha_2)
  $%
   \item
$h^3 \ee^{\frac{1}{2} h A}  \times \left\{r[B(\tfrac12 h),[A,B(\tfrac12 h)]]+\tfrac16 q^3 B^3(\tfrac12 h)\right\}  \times \ee^{\frac{1}{2} h A}u_0\  
\leadsto h^3r[1,-1]\cdot\grad f_2(\alpha_2,\alpha_2)+\frac{1}{3!}h^3v_2^3f_3(\alpha_2,\alpha_2,\alpha_2)$%
  \item  
$\frac12 h^3 p^3\ee^{h\tau A} \times B((1-\tau)h)  \times \ee^{h(1-2\tau)A} \times B^2(\tau h) \times \ee^{h\tau A}u_0\ \leadsto \frac{1}{2!}h^3v_3^2v_1f_3(\alpha_3,\alpha_1,\alpha_1)$%
  \item
  $h^3 p^2q\,\ee^{h\tau A} \times    B((1-\tau)h) \times \ee^{h(\frac12-\tau)A} \times B(\tfrac12 h) \times \ee^{h(\frac12-\tau)A}  \times  B(\tau h) \times \ee^{h\tau A}u_0\ \leadsto h^3v_3v_2v_1f_3(\alpha_3,\alpha_2,\alpha_1)$%
\item  
$\frac12 h^3 p^3\,\ee^{h\tau A} \times  B^2((1-\tau)h)  \times \ee^{h(1-2\tau)A} \times B(\tau h) \times \ee^{h\tau A}u_0\leadsto \frac{1}{2!}h^3v_3^2v_1f_3(\alpha_3,\alpha_3,\alpha_1)$%

\item[] {\bf $\mathcal{O}(h^4)$ terms}

  \item
  $\frac12 h^4 p^2q^2\ee^{h\tau A} \times   B((1-\tau)h) \times \ee^{h(\frac12-\tau)A} \times  B^2(\tfrac12 h) \times \ee^{h(\frac12-\tau)A}  \times B(\tau h) \times \ee^{h\tau A}u_0\  \leadsto \frac{1}{2!}h^4v_3v_2^2v_1f_4(\alpha_3,\alpha_2,\alpha_2,\alpha_1)$  
\item
$ \frac14 h^4 p^2q^2 \ee^{h\frac12 A} \times B^2(\tfrac12 h) \times \ee^{h(\frac12-\tau)A} \times B^2(\tau h) \times \ee^{h\tau A}u_0
\ 
\leadsto \frac{1}{2!\times2!}h^4v_2^2v_1^2f_4(\alpha_2,\alpha_2,\alpha_1,\alpha_1)$
  \item
  $\frac14 h^4 p^2 q^2 \ee^{h\tau A} \times   B^2((1-\tau)h) \times  \ee^{h(\frac12-\tau)A} \times B^2(\tfrac12 h) \times \ee^{h\frac12A}u_0\ 
  \leadsto \frac{1}{2!\times2!}h^4v_3^2v_2^2f_4(\alpha_3,\alpha_3,\alpha_2,\alpha_2)
  $
  \item
  $\frac12 p^3q\,h^4 \ee^{h\tau A} \times    B((1-\tau)h) \times \ee^{h(\frac12-\tau)A} \times B(\tfrac12 h) \times \ee^{h(\frac12-\tau)A}  \times  B^2(\tau h) \times \ee^{h\tau A}u_0\ 
   \leadsto \frac{1}{2!}h^4v_3v_2v_1^2f_4(\alpha_3,\alpha_2,\alpha_1,\alpha_1)
  $
  \item
  $\frac12 h^4 p^3q\,\ee^{h\tau A} \times   B^2((1-\tau)h) \times  \ee^{h(\frac12-\tau)A} \times B(\tfrac12 h) \times \ee^{h(\frac12-\tau)A}  \times B(\tau h) \times \ee^{h\tau A}u_0\ 
  \leadsto \frac{1}{2!}h^4v_3^2v_2v_1f_4(\alpha_3,\alpha_3,\alpha_2,\alpha_1)$
  \item
$ \frac16 h^4 p^3q\,\ee^{h\tau A} \times  B^3((1-\tau)h) \times \ee^{h(\frac12-\tau)A} \times B(\tfrac12 h) \times  \ee^{\frac1 2 h A}u_0\ 
\leadsto \frac{1}{3!}h^4v_3^3v_2f_4(\alpha_3,\alpha_3,\alpha_3,\alpha_2)
$
\item
$ \frac16 h^4 p^3q\,\ee^{h\frac12 A} \times B(\tfrac12 h) \times \ee^{h(\frac12-\tau)A} \times B^3(\tau h) \times \ee^{h\tau A}u_0\ 
\leadsto \frac{1}{3!}v_2v_1^3f_4(\alpha_2,\alpha_1,\alpha_1,\alpha_1)$
\item
$\frac{1}{24} h^4 p^4 \ee^{h\tau A}  \times B^4((1-\tau)h)   \times \ee^{h(1-\tau)A}u_0\ \leadsto \frac{1}{4!}v_3^4f_4(\alpha_3,\alpha_3,\alpha_3,\alpha_3)
$
   \item
 $\frac{1}{24} h^4 p^4 \ee^{h(1-\tau)A} \times B^4(\tau h) \times \ee^{h\tau A}u_0\ 
 \leadsto \frac{1}{4!}h^4v_1^4f_4(\alpha_1,\alpha_1,\alpha_1,\alpha_1)$
 \item  
$\frac16 h^4 p^4 \ee^{h\tau A} \times  B^3((1-\tau)h)  \times \ee^{h(1-2\tau)A} \times B(\tau h) \times \ee^{h\tau A}u_0\ 
\leadsto \frac{1}{3!}h^4v_3^3v_1f_4(\alpha_3,\alpha_3,\alpha_3,\alpha_1)$
\item  
$\frac14 h^4 p^4 \ee^{h\tau A} \times  B^2((1-\tau)h)  \times \ee^{h(1-2\tau)A} \times  B^2(\tau h) \times \ee^{h\tau A}u_0\ 
\leadsto \frac{1}{2!\times2!}h^4v_3^2v_1^2f_4(\alpha_3,\alpha_3,\alpha_1,\alpha_1)
$
\item  
$\frac16 h^4 p^4 \ee^{h\tau A} \times  B((1-\tau)h)  \times \ee^{h(1-2\tau)A} \times B^3(\tau h) \times \ee^{h\tau A}u_0\ 
\leadsto \frac{1}{3!}h^4v_3v_1^3f_4(\alpha_3,\alpha_1,\alpha_1,\alpha_1)
$
    \item
 $ h^4 p\,\ee^{h\tau A} \times  B((1-\tau)h) \times \ee^{h(\frac12-\tau)A} \times  \left\{r[B(\tfrac12 h),[A,B(\tfrac12 h)]]+\tfrac16 q^3 B^3(\tfrac12 h)\right\} \times \ee^{\frac12hA}\\ 
 \leadsto 
 h^4rv_3[0,1,-1]\cdot\grad f_3(\alpha_3,\alpha_2,\alpha_2)+\frac{1}{3!}h^4v_3v_2^3f_4(\alpha_3,\alpha_2,\alpha_2,\alpha_2)
 $
\item
$h^4 p \ee^{h \frac12 A} \times  \left\{r[B(\tfrac12 h),[A,B(\tfrac12 h)]]+\tfrac16 q^3 B^3(\tfrac12 h)\right\} \times \ee^{h(\frac12-\tau)A} \times B(\tau h) \times \ee^{h\tau A}u_0\\ 
 \leadsto 
h^4 r v_1[1,-1,0]\cdot\grad f_3(\alpha_2,\alpha_2,\alpha_1)
 +\frac{1}{3!}h^4v_2^3v_1f_4(\alpha_2,\alpha_2,\alpha_2,\alpha_1)
$
\item
$h^4 \ee^{\frac{1}{2} h A} \times \left\{ \tfrac12 qr B(\tfrac12 h)[B(\tfrac12 h),[A,B(\tfrac12 h)]]+\tfrac12 qr [B(\tfrac12 h),[A,B(\tfrac12 h)]] B(\tfrac12 h)+\tfrac{1}{24} q^4 B^4(\tfrac12 h)\right\}  \times \ee^{\frac{1}{2} h A}u_0\\ 
 \leadsto  \frac{1}{2!}h^4rv_2[1,0,-1]\cdot\grad f_3(\alpha_2,\alpha_2,\alpha_2)+\frac{1}{4!}h^4v_2^4f_4(\alpha_2,\alpha_2,\alpha_2,\alpha_2)$
\end{enumerate}

\end{document}